\documentclass[12pt,a4paper]{elsarticle}
\usepackage{verbatim}
\usepackage{amsmath}
\usepackage{amssymb}
\usepackage{dsfont}
\usepackage{natbib}

\journal{ArXiv}
\biboptions{authoryear}

\begin{document}

\begin{frontmatter}
  \title{On the derivation of the Khmaladze transforms}
  \author{Leigh A Roberts}
    \address{School of Economics and Finance,\\Victoria University of Wellington,\\Wellington, New Zealand}

\begin{abstract}

Some 40 years ago Khmaladze introduced a transform
which greatly facilitated the distribution free goodness of fit testing of statistical hypotheses.
In the last decade, he has published a related transform, 
broadly offering an alternative means to the same end.

The aim of this paper is to derive these transforms using relatively elementary means, making some simplifications, but losing little in the way of generality.  In this way it is hoped to make these transforms more accessible and more widely used in statistical practice.  We also propose a change of name of the second transform to the Khmaladze rotation, in order to better reflect its nature.

\end{abstract}

\begin{keyword}
\texttt{Khmaladze transform\sep distribution free \sep goodness of fit\sep linear algebraic derivation\sep projection\sep reflection}
\end{keyword}
\end{frontmatter}

\tableofcontents

\section{Introduction}

The generic goodness of fit problem lies in separating possible outcomes of a statistical experiment into a finite number of cells, noting expected and observed frequencies for each cell, and testing for smallness of differences between them.  The classic test for this is the chi squared test, which is distribution free, in the sense that the distribution of the chi squared statistic does not depend on the distribution
generating the data,
provided only that the hypothetical distribution be fully specified.

The vital contribution made by Khmaladze some 40 years ago
was to allow distribution free goodness of fit tests for compound hypotheses \citep{khmaladze1979a,khmaladze1981a}.
Using the Khmaladze Transform (hereafter KT, or occasionally the `first' transform, or `KT1'), one could test for whether the data could plausibly arise from a given distributional family.
This more or less corresponds to what is done in practice:
only rarely would a statistician wish to test for a specific parametric value within the distributional family of interest, which was the only possibility available before the introduction of the KT.

Papers utilising the KT are cited in \cite{li2009a}, \cite{koul-swordson2011a} and \cite{kim2016a}, {\em i.a}.  It is however clear that the uptake by the statistical community of such an important conceptual advance in goodness of fit testing has been slow.

Within the last decade, Khmaladze has published another, and ostensibly simpler, transform to test goodness of fit in a distribution free manner \citep{khmaladze2013a,khmaladze2016a}.
The new transform has tentatively been labelled as the second Khmaladze transform, generally referred to below as the second transform or `KT2'.
We suggest that this second transform be relabelled as the `Khmaladze rotation'.

Adoption of the second transform also appears to be slow.  Of the papers
\cite{dumitrescu-khmaladze2019a},
\cite{kennedy2018a},
\cite{khmaladze2017a},
\cite{nguyen-ttm2017b,nguyen-ttm2017a} and
\cite{roberts2019a},
only Kennedy applies the KT2 to real data.
He is also the only one of these authors to suggest a wider use of the second KT to choose an optimal model amongst competing models.

Both transforms involve projections of the empirical process.
The first transform projects onto a Brownian motion (BM), obtained essentially by regressing the incremental empirical distribution function (EDF) on the `future';
what this cryptic description really means is that the goodness of fit test assumes that all the data is to hand, viz.\ that the statistical experiment is complete, when the goodness of fit test is carried out.  The second transform utilises the same types of projections to change from one empirical process to another, without losing any statistical information.  The test of goodness of fit can then be carried out in one statistical framework or the other, whichever is more convenient.

The purpose of this paper is not to apply these transforms to data analysis, but to spread the word about these elegant and potentially very useful transforms, by explaining their genesis more simply and intuitively. The idea is to discretise the empirical processes being tested, so that operators and operands in functional space become matrices and vectors.  We then develop the ideas underlying the transforms in $N-$dimensional space, hoping that the straightforward linear algebraic approach will make the underlying reasoning transparent.  The loss of generality is in fact slight, and such discretisation of the underlying spaces, operators and operands is no more than one would impose for computing purposes.

Our starting point is to decompose the chi squared statistic.  Mimicking the approach taken in \cite{khmaladze2013a}, we project a standard normal vector onto a vector having the covariance structure of the chi squared components; then we consider rotation or reflection from one chi squared statistic to another.  This already gives the two elements underlying the second transform, viz.\ projection and rotation/reflection.  The idea of projection was not new: the fitting of linear regression underpinning the first transform is a projection of asymptotically normal regressors to model the increment of the empirical process.  It is the rotation, or more properly reflection, from one chi squared statistic to another which defines the underlying rationale of the second transform, and differentiates it from the first transform.  Our basic task is to find the projection and rotation operators of the second transform, moving from one (discretised) empirical process to another.

Following the definitions of the projection and reflection operators applied to the chi squared statistic, we then make the ostensibly slight adjustment to operate on the empirical process, essentially the numerator of the chi squared statistic.  Starting with Brownian Motion (BM) in a discrete time $P$, we recast operators as matrices and functional operands as vectors, possibly semi-infinite in length. We project the BM to a Brownian bridge (BB) in time $P$, and rotate from one BB to another, in time $R$ say.  Allowing for parameters to be estimated is allowed for by further projections, most easily seen by changing from point parametrisation to functional parametrisation of the empirical process.

The resulting $q-$projected BMs may be rotated from one empirical process to another, say from time $P$ to time $R$.  We are still working with discrete distributions; but within that limitation, we have proved the validity of the KT2, whereby a goodness of fit test may be effected either in the $P$ space or the $R$ space as convenient, and no statistical information is lost in rotation from one to the other.

A further short section extends the second transform to higher dimensional distributions $P$, and illustrates the outworking for the colour blind problem, for which we largely follow the first part of \cite{dumitrescu-khmaladze2019a}.

Khmaladze
has largely adopted this discretised approach
to elucidate and prove the original transform in the framework of a simple mortality investigation \citep[ch.\ 7]{khmaladze2013b}.
In the final section of this paper we flesh out his development and provide additional comments.

\section{The second Khmaladze transform, or the Khmaladze rotation}

\subsection{The statistical setup\label{s 2}}

Given a random variable X with distribution function $F(x)$,
it is not essential but notationally convenient to suppose that the support of $X$ is bounded away from minus infinity: suppose there exists a finite number $M_F$ such that $X>M_F$ with probability one. Then we
define a grid of points $M_F=x_1<x_2<\ldots<x_N<x_{N+1}=\infty$, and let the $j$th cell be defined by $X\in[x_j,x_{j+1})$ for $1\leq j\leq N$, with associated probabilities $\int_{[x_j,x_{j+1})}dF(x)=p_j$.  Should there be an atom of $X$ at the grid point $x_j$, then the saltus at $x_j$ will be included in $p_j$ but not in $p_{j-1}$.  Noting that $\sum_{j=1}^Np_j=1$, we further define the approximating discrete distribution function $P(x)$, with atoms at $\{x_j\}_{j=1}^N$.

In simpler terms, the non-decreasing step function $P:\mathbb R\to[0,1]$, is piecewise constant apart from steps occurring at points $x_1,x_2,\ldots,x_N$.
The step or saltus at $x_j$ is to be $p_j>0$;
the function $P$ is to be continuous from the right; $P(x)=0$ when $x<x_1$; and $P(x)=1$ when $x\geq x_N$.
The point of supposing the existence of such a lower bound as $M_F$ is merely for the convenience of not regarding $-\infty$ as a possible value for the corresponding putative random variable: it is not an essential step.

In like vein, suppose a random variable $Y$ with distribution function $G(y)$, similarly bounded away from minus infinity, with a grid $M_G=y_1<y_2<\ldots<y_N<y_{N+1}=\infty$, and let the $j$th cell be defined by $Y\in[y_j,y_{j+1})$ for $1\leq j\leq N$, with probabilities $\int_{[y_j,y_{j+1})}dG(y)=r_j$.

The aim of the second Khmaladze transform or Khmaladze rotation is to rotate the stochastic process with distribution function $F$ to that with distribution function $G$.  Our simplified approach is to approximate $F$ and $G$ by the discrete distributions $P$ and $R$ respectively, and rotate from $P$ to $R$.

This latter is more easily visualised than the full rotation from $F$ to $G$, because all operations may be effected by standard linear algebraic procedures.
Resulting vectors and matrices are approximations to the integral operators and functional vectors occurring in the full rotation.
The outworkings below will be set out for finite dimension $N$, but would be essentially unchanged if $N$ were countably infinite, with vectors of semi-infinite length.

\subsection{Rotating from one $\chi^2$ test to another\label{s 1}}

\subsubsection{Projection of chi square statistics}

More or less following \cite{khmaladze2013a}, we normalise the observed frequencies and sum to obtain the conventional chi-squared statistic.
For a sample size $n$ and a total number $N$ of cells, define
\begin{equation}
Y^\prime_j=\frac{\nu_j-np_j}{\sqrt{np_j}}\quad\mbox{for}\ 1\leq j\leq N
\label{e 9}
\end{equation}
with $\nu_j$ the observed frequency in the $j$th cell, and $p_j$ the probability of a data point falling in the $j$th cell.

Assuming that $p_j>0$ for all $j$, we set
$Y^\prime=\begin{pmatrix}Y^\prime_1&Y^\prime_2&\ldots&Y^\prime_N\end{pmatrix}^T$,
$p=\begin{pmatrix}p_1&p_2&\ldots&p_N\end{pmatrix}^T$ and 
$\sqrt p=\begin{pmatrix}\sqrt p_1&\sqrt p_2&\ldots&\sqrt p_N\end{pmatrix}^T$.
Note that $\sqrt p^T\sqrt p=\sum_{j=1}^n p_j=1$.
Summing the squares of the statistics in \eqref{e 9} produces the conventional chi-squared statistic
\begin{equation}
\chi^2=\sum_{j=1}^n\frac{(O_j-E_j)^2}{E_j}=\sum_{j=1}^n{Y^\prime_j}^2={Y^\prime}^TY^{\prime}
\label{e 8}
\end{equation}
\cite{chibisov1971a}, for example, discusses the distribution of the chi-squared statistic when cell boundaries are determined in light of the data; but we assume boundaries to be fixed independently of the data.

The frequencies $\nu_j$ have a multinomial distribution, with 
covariance structure given by $\mbox{Var}(\nu_j)=np_j(1-p_j)$, $\mbox{Cov}(\nu_j,\nu_k)=-np_jp_k$ for $j\neq k$.
The covariance matrix of $Y^\prime$ reduces to
\begin{equation}
E\ Y^\prime {Y^\prime}^T
=\begin{pmatrix}1-p_1&-\sqrt p_1\sqrt p_2&\ldots&-\sqrt p_1\sqrt p_N\\
-\sqrt p_2\sqrt p_1&1-p_2&\ldots&-\sqrt p_2\sqrt p_N\\
&&\ldots&\\
-\sqrt p_N\sqrt p_1&-\sqrt p_N\sqrt p_2&\ldots&1-p_N\\  \end{pmatrix}
=I-\sqrt p\sqrt p^T
\label{e 1}
\end{equation}

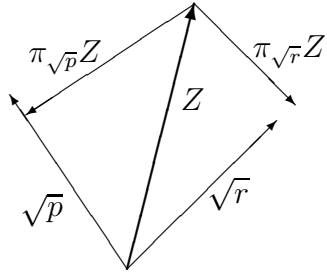
\begin{figure}
\begin{center}
\setlength{\unitlength}{1pt}
\begin{picture}(120,100)
\thicklines
\put(50,0){\vector(1,4){25}}
\thinlines
\put(50,0){\vector(1,1){56}}
\put(75,100){\vector(1,-1){38}}
\put(50,0){\vector(-2,3){44}}
\put(75,100){\vector(-3,-2){63}}
\put(80,25){$\sqrt r$}
\put(10,20){$\sqrt p$}
\put(70,60){$Z$}
\put(97,80){$\pi_{\sqrt r}Z$}
\put(13,78){$\pi_{\sqrt p}Z$}
\end{picture}
\caption{Projections of $Z$ perpendicular to the unit vectors $\sqrt p$ and $\sqrt r$\label{f 1}}
\end{center}
\end{figure}

Define $\pi_a(b)$ to be a projection from $b$ perpendicular to $a$, as illustrated in Figure \ref{f 1}. First note that
\begin{equation}
\pi_{\sqrt p}=I-\sqrt p\sqrt p^T
\qquad\qquad\pi_{\sqrt p}^T=\pi_{\sqrt p}\qquad\qquad\pi_{\sqrt p}^2=\pi_{\sqrt p}
\label{e 32}
\end{equation}
Setting $Z\sim{\cal N}(0,I)$, define the Gaussian vector $Y$
\[
Y=Z-\sqrt p{\sqrt p}^TZ=\pi_{\sqrt p}Z\qquad\qquad \sqrt p^T\,Y=0\qquad\qquad Y^T=Z^T\pi_{\sqrt p}
\]
It is well known that the chi-squared statistic in \eqref{e 8} has the $\chi^2_{n-1}$ limiting distribution, with mean $n-1$.  For completeness, and for comparison with later work, we note that the asymptotic limiting statistic has the same mean:
\[
E\ Y^TY=E\ Z^T\pi_{\sqrt p}\pi_{\sqrt p}Z=E\ Z^T\pi_{\sqrt p}Z=E\ Z^T\left(I-\sqrt p\sqrt p^T\right)Z
\]
\begin{equation}
=E\ Z^TZ-E\ Z^T\sqrt p\sqrt p^TZ=E\ Z^TZ-E\ \sqrt p^TZZ^T\sqrt p
=E\ Z^TZ-\sqrt p^T\sqrt p=n-1
\label{e 14}
\end{equation}
We are however more concerned with the covariance of $Y$:
\begin{equation}
E\ YY^T=E\ \pi_{\sqrt p}ZZ^T\pi_{\sqrt p}=\pi_{\sqrt p}E(ZZ^T)\pi_{\sqrt p}=\pi_{\sqrt p}^2=\pi_{\sqrt p}
\label{e 16}
\end{equation}
or
\[
E\ YY^T=
\begin{pmatrix}1&0&\ldots&0\\0&1&\ldots&0\\&&\ldots&\\0&0&\ldots&1\end{pmatrix}
-\begin{pmatrix}\sqrt p_1\\\sqrt p_2\\\ldots\\\sqrt p_N\end{pmatrix}
\begin{pmatrix}\sqrt p_1&\sqrt p_2&\ldots&\sqrt p_N\end{pmatrix}
\]
\[
=\begin{pmatrix}1-p_1&-\sqrt p_1\sqrt p_2&\ldots&-\sqrt p_1\sqrt p_N\\
-\sqrt p_2\sqrt p_1&1-p_2&\ldots&-\sqrt p_2\sqrt p_N\\
&&\ldots&\\
-\sqrt p_N\sqrt p_1&-\sqrt p_N\sqrt p_2&\ldots&1-p_N\end{pmatrix}
\]
in agreement with the covariance of $Y^\prime$ in \eqref{e 1}.  As $n\to\infty$, $Y^\prime$ tends weakly to $Y$, i.e.,
the (cumulative) distribution function of $Y^\prime$ tends to that of $Y$, at all points of continuity of the latter.  In simpler terms,
we may asymptotically approximate $Y^\prime$ by the normally distributed $Y$. This is certainly so if the infimum of $\{p_j:1\leq j\leq N\}$ is bounded away from zero.

\subsubsection{From a projection operator to a reflection operator}

Given vectors of unit length $\sqrt p$ and $\sqrt r$, the reflection operator $U_{\sqrt p,\sqrt r}=U_0$ is defined as
\[
U_{\sqrt p,\sqrt r}=U_0=I-c_0(\sqrt p-\sqrt r)(\sqrt p-\sqrt r)^T
\]
in which the constant $c_0$ is given by
\[
c_0=\frac2{\parallel\sqrt p-\sqrt r\parallel^2}=\frac1{1-<\sqrt p,\sqrt r>}
\]
and in turn the
norm is given by $\|s\|^2=\ <s,s>\ =s^Ts$ and the
inner product as $<s,t>\ =s^Tt$.
The reflection operator swaps $\sqrt p$ and $\sqrt r$ around, while leaving unchanged all vectors orthogonal to both of them.
We note that $U_0^T=U_0$, $U_0^2=I$ and $U_0\pi_{\sqrt p}U_0=\pi_{\sqrt r}$.

The covariance of $U_0Y$ is
\[
\mbox{cov}(U_0Y)=E\,U_0Y(U_0Y)^T=U_0\,E\,YY^T\,U_0=U_0\pi_{\sqrt p}U_0=\pi_{\sqrt r}
=I-\sqrt r\sqrt r^T
\]

Consider a second stochastic vector of length $N$ defined by $T=\pi_{\sqrt r}Z_1$ where $Z_1\sim{\cal N}(0,I)$ and $Z,Z_1$ are independent.

We may consider $T^{\,\prime}$ to be defined from a normalised multinomial variate as in \eqref{e 9}, with $p_j$ replaced by $r_j$; again we may asymptotically approximate $T^{\,\prime}$ by the normally distributed $T$.
Thus $T$ and $U_0Y$ have the same distribution, since they are Gaussian with identical means and covariance.

So given two stochastic processes and a common number $N$ of cells with cell probabilities $p$ and $r$, we have a means of rotating from one process to another, or at least from one chi-squared test to another.  The projections of $Z$ and $Z_1$ to give $Y$ and $T$ respectively are not reversible, and so lose information.  In contrast, the transform from $Y$ to $T$ (we are working asymptotically, and assume normality) is reversible.  The goodness of fit test may be more conveniently carried out for one process than the other, and there is no loss of statistical information in rotating from one stochastic process to another.

\subsection{From the chi-squared statistic to the empirical process}

We have normalised the multinomial numerator in \eqref{e 9} in order to define the rotation $U_0$ between the stochastic processes $Y^\prime$ and $T^{\,\prime}$, or rather between their asyptotic limits $Y$ and $T$.  Empirical processes in practice assume the form of observed minus expected frequencies, and our first task is to remove the standard deviation in the denominator from \eqref{e 9}, although we retain the factor of $\sqrt n$ needed for sensible scaling of the empirical process.

Accordingly we define $N\times N$ diagonal matrices $D_p$ containing $\{p_j\}_1^N$, in the given order, along the diagonal: that is, $D_p=\mbox{diag}(p)$.
The square root of $D_p$ is unambiguously defined, with elements $\sqrt {p_j}$ down the diagonal, denoted either by $D_{\sqrt p}$ or $D_p^{1/2}$.
Thus the vector $D_p^{1/2}Y^\prime$ contains elements $(\nu_j-np_j)/\sqrt n$, for $j=1,\ldots,n$, although we generally work instead with the limiting Gaussian vector $D_p^{1/2}Y$.
The diagonal $N\times N$ matrix $D_r$ is defined analogously for probabilities $r$.

Define a lower triangular $N\times N$ matrix $J$ by setting all elements on and below the diagonal to unity, and the elements above the diagonal to zero.  This `accumulating' matrix has the effect of cumulating a column vector: the $j$th element of $Y$, for instance, is $Y_j$, while
the $j$th element of $JY$ is $\sum_{k=1}^jY_k$.
Further, let $j_k^T$ denote the $k$th row of the accumulation matrix $J$, so that the column vector $j_k$ consists of $k$ unities followed by $N-k$ zeroes.

A caution is in order here.  A common convention is to denote
random variables by capital letters, with possible or realised (sample) values denoted by the corresponding small letter.  We do not necessarily keep to this convention here.

\subsubsection{Brownian Motion in time $P(x)$}

Consider again the vector $Z\sim{\cal N}(0,I)$.
The covariance matrix of $D_p^{1/2}Z$ is given by
\[
\mbox{Cov}\left(D_{\sqrt p}Z\right)=E\,D_{\sqrt p}ZZ^TD_{\sqrt p}=D_p
\]
The elements of the vector $D_p^{1/2}Z$ are increments of Brownian motion (BM) in the time $P(x)$.  Integrating or summing that process yields BM in time $P(x)$.
The stochastic vector $JD_p^{1/2}Z$ has the distribution of a BM with respect to the time $P(x)$, or in the time $P(x)$; and its covariance matrix is $J\,D_p\,J^T$.

We define $\vec{\Delta w_P}=D_p^{1/2}Z$,
the vectorised increments of BM in time $P$; we further set $\vec{w_P}=JD_p^{1/2}Z$, the vectorised BM in time $P$, or the vectorised $P$ BM.

Suppose $x\in[x_k,x_{k+1})$.  The conventional expression of BM in time $P$, evaluated at time $x$, would then be
\begin{equation}
w_P(x)
=\sum_{j\leq k}\sqrt{p_j}\,Z_j=j_k^T\,D_{p}^{1/2}\,Z
\label{e 41}
\end{equation}
which is normally distributed, with mean zero and variance $\int_{(-\infty,x]}dP(y)=\sum_{j\leq k}p_j$.

Further assume that $x^{\,\prime}\in[x_l,x_{l+1})$, where $k\leq l$.
The covariance of $w_P(x)$ and $w_P(x^{\,\prime})$
is given by
\[
  \mbox{Cov}\left(w_P(x),w_P(x^{\,\prime})\right)
  =\mbox{Cov}\left(j_k^T\,D_{p}^{1/2}\,Z,j_l^T\,D_{p}^{1/2}\,Z\right)
\]
\[
=j_k^T\,D_{p}^{1/2}\,EZZ^T\,D_{p}^{1/2}\,j_l=j_k^T\,D_{p}\,j_l=\sum_{j\leq k}p_j
=\int_{-\infty}^{{\scriptsize\mbox{min}}(x,x^{\,\prime})}dP(y)
=P\left(\mbox{min}(x,x^{\,\prime})\right)
\]
which is the standard expression for the covariance of BM
in time $P(x)$ for any distribution function $P(x)$.

Now suppose that $x\in[x_k,x_{k+1})$ and $x+\Delta x\in[x_m,x_{m+1})$, where $k\leq m$.
Conventional increments of $w_P(x)$ are given by
\[
\Delta w_P(x)=w_P(x+\Delta x)-w_P(x)=\sum_{k<j\leq m}\sqrt{p_j}\,Z_j=(j_m-j_k)^T\,D_{p}^{1/2}\,Z
\]
where $\Delta w_P(x)$ is again normally distributed, with mean zero and variance $\int_x^{x+\Delta x}dP(y)=\sum_{k<j\leq m}p_j$.

\subsubsection{The Brownian Bridge in time $P(x)$}

In the last section,
we rescaled the standard normal column vector $Z$ and interpreted the elements of $D_p^{1/2}Z$ as increments of BM in time $P(x)$.
In like manner, the column vector $D_p^{1/2}Y$ contains
increments of the Brownian Bridge (BB) process $v_P(x)$ in time $P(x)$, so that the BB process in time $P(x)$ is $JD_p^{1/2}Y$.

As previously,
we define $\vec{\Delta v_P}=D_p^{1/2}Y$,
the vectorised increments of the BB in time $P$; and $\vec{v_P}=JD_p^{1/2}Y$, the vectorised BB in time $P$.

By analogy with \eqref{e 41} above, and still assuming that $x_k\leq x<x_{k+1}$, we have that
\[
  v_P(x)
=\sum_{j\leq k}\sqrt{p_j}\,Y_j=j_k^T\,D_{p}^{1/2}\,Y
\]
which is normally distributed, with mean zero; but the covariance structure of the $P$ BB is now more complicated.

Again supposing that  $x^{\,\prime}\in[x_l,x_{l+1})$, with $k\leq l$,
the covariance of $v_P(x)$ and $v_P(x^{\,\prime})$
is given by
\begin{equation}
  \mbox{Cov}\left(v_P(x),v_P(x^{\,\prime})\right)
  =\mbox{Cov}\left(j_k^T\,D_{p}^{1/2}\,Y,j_l^T\,D_{p}^{1/2}\,Y\right)
\label{e 44}
\end{equation}
\[
=j_k^T\,D_{p}^{1/2}\,EYY^T\,D_{p}^{1/2}\,j_l
=j_k^T\,D_{p}^{1/2}\left(I-\sqrt p\sqrt p^T\right)D_{p}^{1/2}\,j_l
\]
\[
=j_k^T\,D_{p}\,j_l-j_k^T\,pp^T\,j_l
=\sum_{j\leq k}p_j-\sum_{j\leq k}p_j\sum_{j\leq l}p_j
\]
\[
=\int_{-\infty}^{{\scriptsize\mbox{min}}(x,x^{\,\prime})}dP(y)-\int_{-\infty}^xdP(y)\int_{-\infty}^{x^{\,\prime}}dP(y)
=P\left(\mbox{min}(x,x^{\,\prime})\right)-P\left(x\right)P\left(x^{\,\prime}\right)
\]
which is the standard expression for the covariance of a BB
in time $P(x)$, or a $P$ BB for short, for any distribution function $P(x)$.

When
$x+\Delta x\in[x_m,x_{m+1})$
the increments of $v_P(x)$ are given by
\[
\Delta v_P(x)=v_P(x+\Delta x)-v_P(x)=\sum_{k<j\leq m}\sqrt{p_j}\,Y_j=(j_m-j_k)^T\,D_{p}^{1/2}\,Y
\]
where $\Delta v_P(x)$ is again normally distributed, with mean zero.
The variance of $\Delta v_P(x)$ is
\[
\mbox{Var}(\Delta v_P(x))=(j_m-j_k)^T\,D_{p}^{1/2}\,EYY^T\,D_{p}^{1/2}(j_m-j_k)
\]
\[
=(j_m-j_k)^T\,D_{p}^{1/2}\left(I-\sqrt p\sqrt p^T\right)D_{p}^{1/2}(j_m-j_k)
\]
\[
=(j_m-j_k)^T\,D_{p}(j_m-j_k)-(j_m-j_k)^T\,pp^T\,(j_m-j_k)
\]
\begin{equation}
=\sum_{k<j\leq m}p_j-\left(\sum_{k<j\leq m}p_j\right)^2
\label{e 49}
\end{equation}

\subsubsection{Projection from BM to BB in time $P(x)$}

Let $q_0$ be a column vector of length $N$ in which every element is unity.
\[
Y=\pi_{\sqrt p}Z=Z-\sqrt p{\sqrt p}^TZ
=\left(I-D_p^{1/2}q_0q_0^TD_p^{1/2}\right)Z
\]
\[
D_p^{1/2}Y
=\left(I-D_pq_0q_0^T\right)D_p^{1/2}Z
=\Pi_P^{q_0}D_p^{1/2}Z
\]
\[
\vec{\Delta v_P}
=\left(I-D_pq_0q_0^T\right)\vec{\Delta w_P}
=\Pi_P^{q_0}\,\vec{\Delta w_P}
\]
Working in the primal space with $v_P$ and $w_P$, the projection operator $\Pi_P^{q_0}$ projects $P$ BM $\vec{w_P}$ onto $q_0-$projected $P$ BM, which is just the $P$ BB $\vec{v_P}$ \citep{khmaladze2016a}.  Abbreviating for the moment by setting $\Pi_P^{q_0}=\Pi$ (to be generalised later in \eqref{e 22} on p.\ \pageref{e 22}), we have
\begin{equation}
\Pi_P^{q_0}=\Pi=I-D_pq_0q_0^T\qquad\Pi^2=\Pi\qquad\Pi D_p\Pi^T=D_p-D_pq_0q_0^TD_p=\Pi D_p=D_p\Pi^T
\label{e 21}
\end{equation}
Noting that $q_0^TD_pq_0=1$,
proofs are as follows
\[
\Pi^2=\left(I-D_pq_0q_0^T\right)\left(I-D_pq_0q_0^T\right)=I-2D_pq_0q_0^T+D_pq_0q_0^TD_pq_0q_0^T=I-D_pq_0q_0^T=\Pi
\]
\[
\Pi D_p\Pi^T=\left(I-D_pq_0q_0^T\right)D_p\left(I-q_0q_0^TD_p\right)
\]
\[
=D_p-2D_pq_0q_0^TD_p+D_pq_0q_0^TD_pq_0q_0^TD_p=D_p-D_pq_0q_0^TD_p
\]

\subsection{From point parametric to function parametric form for the stochastic processes}

Now we change to function parametric form for the BM and BB.

Let $\Phi=\{\phi\}$ be a family of functions $\phi(x)\in L^2_P$, i.e.\ functions $\phi$ such that $\int\phi(x)^2dP(x)<\infty$.
We set $\vec\phi=\begin{pmatrix}\phi_1&\phi_2&\ldots&\phi_N\end{pmatrix}^T$, where $\phi_j=\phi(x_j)$.
In an abuse of notation we shall often write $\phi$ for $\vec\phi$, since there seems little likelihood of confusing the function $\phi(x)$ and the vector of its non-zero values at the atoms of $P$.
The finiteness condition for the norm of members of $\Phi$ reduces to $\|\phi\|^2_P=\int\phi(x)^2dP(x)=\vec\phi^{\ T}\,D_p\,\vec\phi=\phi^{\ T}\,D_p\,\phi<\infty$, which ceases to be vacuous if $N$ is allowed to be countably infinite, given that the values $\phi_j$ are finite.

Similarly we define a family $\Psi$ of functions $\psi(y)\in L^2_R$, with vector of values $\vec\psi=\begin{pmatrix}\psi_1&\psi_2&\ldots&\psi_N\end{pmatrix}^T$, where $\psi_j=\psi(y_j)$.
Again writing $\psi$ for $\vec\psi$, we are imposing the analogous constraint on the norm, viz.\ $\|\psi\|^2_R=\int\psi(x)^2dR(x)=\vec\psi^{\ T}\,D_R\,\vec\psi=\psi^{\ T}\,D_R\,\psi<\infty$.

\subsubsection{Projection operators in the primal and dual spaces}

Writing $\vec\phi^{\ T}\vec{\Delta v_P}=\phi^{\ T}\vec{\Delta v_P}$ then, we have
\[
\phi^{\ T}\vec{\Delta v_P}=\phi^{\ T}\Pi_P^{q_0}\vec{\Delta w_P}
  =\phi^{\ T}\vec{\Delta w_P}-\phi^{\ T}D_pq_0q_0^T\vec{\Delta w_P}
\]
\begin{equation}
  =\phi^{\ T}\Pi\,\vec{\Delta w_P}=\left(\Pi^T\phi\right)^T\vec{\Delta w_P}
\label{e 10}
\end{equation}
so that $\Pi=I-D_pq_0q_0^T$ is the projection operator in the primal space, acting on $\vec{\Delta w_P}$ to produce $\vec{\Delta v_P}$; and $\Pi^T=I-q_0q_0^TD_p$ is the projection operator in the dual space, acting on the functions $\phi\in\Phi$.
Equation \eqref{e 10} may be rewritten in more conventional fashion as
\[
\int_{-\infty}^{\infty}\phi(x)dv_P(x)=\int_{-\infty}^{\infty}\phi(x)dw_P(x)-\int_{-\infty}^{\infty}\phi(x)q_0(x)dP(x)
\int_{-\infty}^{\infty} q_0(x)dw_P(x)
\]
or more succinctly as
\begin{equation}
v_P^{q_0}(\phi)=v_P(\phi)=w_P(\phi)\ -<\phi,q_0>_P\ w_P(q_0)
\label{e 19}
\end{equation}
in which $v_P(\phi)$ is $q_0$-projected $P$ BM \citep{khmaladze2016a}, or $q_0$-projected BM in time $P(x)$; and where
\[
\int_{-\infty}^{\infty}\phi(x)q_0(x)dP(x)=\int_{-\infty}^{\infty}\phi(x)dP(x)=\ <\phi,q_0>_P
\]
To retrieve the point parametric version $v_P(x)$, set $\phi(s)=\phi_t(s)=\mathds1_{\{s< t\}}$ = the Heaviside function, starting at 1 and dropping to zero at $t$. 
From \eqref{e 19}
\[
\int_{-\infty}^xdv_P(y)=\int_{-\infty}^xdw_P(y)-
\int_{-\infty}^xdP(y)
\int_{-\infty}^{\infty} q_0(x)dw_P(x)
\]
\[
v_P^{q_0}(x)=v_P(x)=w_P(x)-P(x)w_P(\infty)
\]

\subsubsection{Covariance of $v^q_P(\phi)$ and $v^q_P(\widetilde\phi)$ for $q=q_0$}

\[
\mbox{Cov}\left(v_P(\phi),v_P(\widetilde\phi)\right)
=\mbox{Cov}\left(\phi^{\ T}\,\vec{\Delta v_P},{{\widetilde\phi}}\,^{\,T}\,\vec{\Delta v_P}\right)
=\mbox{Cov}\left(\phi^{\ T}\Pi\,\vec{\Delta w_P},{{\widetilde\phi}}\,^{\,T}\Pi\,\vec{\Delta w_P}\right)
\]
\[
=E\ \phi^{\ T}\Pi\,\vec{\Delta w_P}\,\vec{\Delta w_P}^{\,T}\Pi^T{{\widetilde\phi}}=\phi^{\ T}\Pi D_p\Pi^T{{\widetilde\phi}}
\]
\begin{equation}
\mbox{Cov}\left(v_P(\phi),v_P(\widetilde\phi)\right)
=\phi^{\ T}\left(D_p-D_pq_0q_0^TD_p\right){{\widetilde\phi}}
\label{e 11}
\end{equation}
More conventionally perhaps, this would be expressed as
\begin{equation}
Cov\left(v_P(\phi),v_P(\widetilde\phi)\right)
=\int\phi\,\widetilde\phi\, dP-\left(\int\phi\,dP\right)\left(\int\widetilde\phi\, dP\right)
\label{e 25}
\end{equation}
Analogously for the BM with functional parametrisation
\[
\mbox{Cov}\left(w_P(\phi),w_P(\widetilde\phi)\right)
=\mbox{Cov}\left(\phi^{\ T}\,\vec{\Delta w_P},{{\widetilde\phi}}\,^{\,T}\,\vec{\Delta w_P}\right)
\]
\[
=E\ \phi^{\ T}\,\vec{\Delta w_P}\,\vec{\Delta w_P}^T{{\widetilde\phi}}=\phi^{\ T} D_p\,{{\widetilde\phi}}
=\int\phi\,\widetilde\phi\, dP
\]

\subsection{Rotation/reflection operator in the functional (dual) space}

For rotation from one stochastic process to another, we operate on functions in the dual space.  For $\xi$ and $\eta$ functions in $\Phi$ of unit $P$-norm, i.e.\ $\int \xi^2dP=\ <\xi,\xi>_P\ =\|\xi\|^2=\vec\xi^{\ T}D_P\,\vec\xi=\xi^{\ T}D_P\,\xi=1$, and analogously for $\eta$, define the involution
\[
U_{\xi,\eta}=I-c(\xi-\eta)(\xi-\eta)^TD_p
\]
in which
\[
c=\frac2{\|\xi-\eta\|_P^2}=\frac1{1-<\xi,\eta>_P}
\]
and in turn $<\xi,\eta>_P\ =\xi^{\ T}D_P\,\eta$.
Then $U_{\xi,\eta}$ swaps $\xi$ and $\eta$ around, and leaves unchanged any vector orthogonal (or rather $P-$orthogonal) to $\xi$ and $\eta$.  That is,
\[
<\xi,\zeta>_P\ =\ <\eta,\zeta>_P\ =0\quad\Rightarrow\quad U_{\xi,\eta}\zeta=\zeta
\]
The corresponding rotation operator in the primal space is $U_{\xi,\eta}^T$.

Also $U_{\xi,\eta}$ preserves the $P-$norm, i.e., $U_{\xi,\eta}^TD_pU_{\xi,\eta}=D_p$:
for any $N-$vector $\zeta$, $\|U_{\xi,\eta}\zeta\|^2=\zeta^TU_{\xi,\eta}^TD_pU_{\xi,\eta}\zeta=\zeta^TD_p\zeta=\|\zeta\|^2$.

\subsubsection{Rotation from one stochastic process to another when parameters are known}

We set $s_0=q_0$, intending $q_0$ to be associated with the source $P$ distribution, and
$s_0$ to be associated with the target $R$ distribution.  When we allow parameters to be estimated, the vectors $q_j$ and $s_j$ for $0<j\leq K$ will become the (normalised) score functions ($N-$vectors) for the $K$ dimensional parameter $\theta$ in the respective spaces.  For the moment the parameter $\theta$ is assumed known, requiring the use of $q_0$ and $s_0$ alone to define the rotation/reflection.  In the more general case, we shall require a sequence of reflections to define the desired rotation, as set out in \S\ref{s 3} on p.\ \pageref{s 3}.

Define $L=D_r^{1/2}\,D_p^{-1/2}$, so that $LD_pL=D_r$.
Note that for $\psi\in L^2_R$, $L\vec\psi\in L^2_P$: i.e., writing $L\psi$ for $L\vec\psi$,
\[
\left(L\psi\right)^TD_pL\psi=\psi^{\ T}LD_pL\psi=\psi^{\ T}D_r\psi<\infty
\]
Set $U=U_{q_0,Ls_0}$, so that $q_0=ULs_0$.
Then the rotation from one stochastic process to another is given by
\begin{equation}
v^{s_0}_R(\psi\,)=v^{q_0}_P(UL\psi)
\label{e 54}
\end{equation}
Proof is by showing commonality of covariances.
From \eqref{e 11},
\[
Cov\left(v^{q_0}_P(UL\psi),v^{q_0}_P(UL\widetilde\psi)\right)=\psi^TL^TU^T\left(D_p-D_pq_0q_0^TD_p\right)UL\widetilde\psi
\]
\[
=\psi^TL^TD_pL\widetilde\psi-\psi^TL^TU^TD_pULs_0\,s_0^TL^TU^TD_pUL\widetilde\psi
\]
\[
=\psi^TL^TD_pL\widetilde\psi-\psi^TL^TD_pLs_0\,s_0^TL^TD_pL\widetilde\psi
\]
\[
=\psi^TD_r\widetilde\psi-\psi^TD_rs_0\,s_0^TD_r\widetilde\psi
=\psi^T\left(D_r-D_rs_0s_0^TD_r\right)\widetilde\psi
=Cov\left(v^{s_0}_R(\psi),v^{s_0}_R(\widetilde\psi)\right)
\]

\subsection{Estimating parameters}

\subsubsection{The score function}

Recalling the definition of $Y^\prime$ from \eqref{e 9} on p.\ \pageref{e 9}, we define an analogous statistic $\widehat Y_j^\prime$ when the parameter is estimated:
\[
Y_j^\prime=\frac{\nu_j-np_j}{\sqrt{np_j}}\qquad\qquad
\widehat Y_j^\prime=\frac{\nu_j-n\widehat p_j}{\sqrt{n\widehat p_j}}
\]
Estimating $K$ unknown parameters by minimising chi squared \citep[\S30.3]{cramer1946a}, 
\begin{equation}
\widehat Y^\prime=Y^\prime-B(B^TB)^{-1}B^TY^\prime+o_P(1)
\label{e 4}
\end{equation}
$B$ is an $N\times K$ matrix, with rows $j=1,2,\ldots,N$ and columns $k=1,2,\ldots,K$.
\[
B=(B_{jk})=\left(\frac1{\sqrt p_j}\frac{\partial p_j}{\partial\theta_k}\right)
=D_p^{-1/2}\left(\frac{\partial p_j}{\partial\theta_k}\right)
\]
\begin{equation}
=D_p^{1/2}\left(\frac{\partial p_j\big/\partial\theta_k}{p_j}\right)
=D_p^{1/2}\begin{pmatrix}Q_1&Q_2&\ldots&Q_K\end{pmatrix}
\label{e 47}
\end{equation}
in which $Q_k$ is the column $N-$vector with elements $\{\frac{\partial p_j/\partial\theta_k}{p_j}\}_{j=1}^N$ for $k=1,\ldots,K$; in other words, $\begin{pmatrix}Q_1&Q_2&\ldots&Q_K\end{pmatrix}$ is the (non-normalised) score function.
To normalise this function, we note that
\begin{equation}
B^TB=\begin{pmatrix}Q_1^T\\Q_2^T\\\ldots\\Q_K^T\end{pmatrix}D_P\begin{pmatrix}Q_1&Q_2&\ldots&Q_K\end{pmatrix}=\Gamma
\label{e 56}
\end{equation}
where $\Gamma$ is the information matrix.
We may now define normalised score functions as
\[
\begin{pmatrix}q_1&q_2&\ldots&q_K\end{pmatrix}=\begin{pmatrix}Q_1&Q_2&\ldots&Q_K\end{pmatrix}\Gamma^{-1/2}
\]
From \eqref{e 56} then we have
\[
\Gamma^{-1/2}\begin{pmatrix}Q_1^T\\Q_2^T\\\ldots\\Q_K^T\end{pmatrix}D_p\begin{pmatrix}Q_1&Q_2&\ldots&Q_K\end{pmatrix}\Gamma^{-1/2}=I
\]
\[
\begin{pmatrix}q_1^T\\q_2^T\\\ldots\\q_K^T\end{pmatrix}D_p\begin{pmatrix}q_1&q_2&\ldots&q_K\end{pmatrix}=I
\]
So, with $\delta_{jk}$ denoting the Kronecker delta,
\[
q_j^TD_pq_k=\delta_{jk}\qquad\mbox{for}\qquad1\leq j\leq K,1\leq k\leq K
\]
but in fact more is true:
letting the heavy dot denote $\frac{\partial}{\partial\theta_k}$, and noting that $\sum p_j=1$, yields
\[
\sum\overset{\bullet}p_j=0=\sum\frac{\overset{\bullet}p_j}{p_j}\,p_j=Q_k^TD_p\,q_0=q_0^TD_pQ_k
\]
This is true for every $k$, so
\[
q_0^TD_p\begin{pmatrix}Q_1&Q_2&\ldots&Q_K\end{pmatrix}\Gamma^{-1/2}
=q_0^TD_p\begin{pmatrix}q_1&q_2&\ldots&q_K\end{pmatrix}=0
\]
Finally then we may write
\begin{equation}
q_j^TD_pq_k=\delta_{jk}\qquad\mbox{for}\qquad0\leq j\leq K,0\leq k\leq K
\label{e 15}
\end{equation}

\subsubsection{Reassembling the jigsaw}

Fleshing out \eqref{e 4} we find
\[
B(B^TB)^{-1}B^T=D_p^{1/2}\begin{pmatrix}Q_1&Q_2&\ldots&Q_K\end{pmatrix}\Gamma^{-1}\begin{pmatrix}Q_1^T\\Q_2^T\\\ldots\\Q_K^T\end{pmatrix}D_p^{1/2}
\]
\[
=D_p^{1/2}\begin{pmatrix}q_1&q_2&\ldots&q_K\end{pmatrix}\begin{pmatrix}q_1^T\\q_2^T\\\ldots\\q_K^T\end{pmatrix}D_p^{1/2}
\]
\begin{equation}
=\sum_{k=1}^KD_p^{1/2}q_kq_k^TD_p^{1/2}
\label{e 3}
\end{equation}
\begin{equation}
Y=Z-\sqrt p{\sqrt p}^TZ
=\left(I-D_p^{1/2}q_0q_0^TD_p^{1/2}\right)Z=\pi_{\sqrt p}\,Z
\label{e 5}
\end{equation}
From \eqref{e 4} we define the Gaussian $N-$vector
\[
\widehat Y=Y-B(B^TB)^{-1}B^TY
\]
and from \eqref{e 5}
\begin{equation}
\widehat Y=Z-\sqrt p{\sqrt p}^TZ-B(B^TB)^{-1}B^TZ-B(B^TB)^{-1}B^T\sqrt p{\sqrt p}^TZ+o_P(1)
\label{e 2}
\end{equation}
but
\[
B^T\sqrt p=\begin{pmatrix}Q_1^T\\Q_2^T\\\ldots\\Q_K^T\end{pmatrix}D_p^{1/2}\sqrt p=\begin{pmatrix}Q_1^T\\Q_2^T\\\ldots\\Q_K^T\end{pmatrix}D_pq_0=\begin{pmatrix}0\\0\\\ldots\\0\end{pmatrix}
\]
so, from \eqref{e 2},
\begin{equation}
\widehat Y=Z-\sqrt p{\sqrt p}^TZ-B(B^TB)^{-1}B^TZ+o_P(1)
\label{e 6}
\end{equation}
Applying \eqref{e 3} and \eqref{e 5}, and disregarding the residual in \eqref{e 6},
\begin{equation}
  \widehat Y=\left[I-D_p^{1/2}\sum_{k=0}^Kq_kq_k^TD_p^{1/2}\right]Z
  =\widehat\pi_{\sqrt p}Z
\label{e 20}
\end{equation}
Again, in analogy with \eqref{e 32} and \eqref{e 16} on pp.\ \pageref{e 32} and \pageref{e 16}, and utilising \eqref{e 15},
\[
\widehat\pi_{\sqrt p}^T=\widehat\pi_{\sqrt p}\qquad\widehat\pi_{\sqrt p}^2=\widehat\pi_{\sqrt p}\qquad E\ \widehat Y\widehat Y^T=\widehat\pi_{\sqrt p}
\]
From \eqref{e 20} we have that
\[
D_p^{1/2}\widehat Y
=
\left[I-D_p\sum_{k=0}^Kq_kq_k^T\right]
D_p^{1/2}Z
\]
Eschewing the clumsy notation $\vec{\Delta\widehat v_P}$ in favour of the simpler $\Delta\widehat v_P$, we
rewrite the last equation as
\[
  \Delta\widehat v_P
  =D_p^{1/2}\widehat Y
=\left[I-D_p\sum_{k=0}^Kq_kq_k^T\right]
\vec{\Delta w_P}
=\Pi_P^{q}\,\vec{\Delta w_P}
\]
in which $\Pi_P^{q}$ is defined as shown, extending the previously defined $\Pi_P^{q_0}$ in \eqref{e 21} on p.\ \pageref{e 21}, and where $q=\left(q_0,q_1,\ldots,q_K\right)$.

\[
Cov(\Delta\widehat v_P)=cov(D_p^{1/2}\widehat Y)=D_p^{1/2}\widehat\pi_{\sqrt p} D_p^{1/2}
\]

Now we use the abbreviation $\Pi$ for $\Pi_P^{q}$ rather than $\Pi_P^{q_0}$.
Proofs of the relations in \eqref{e 22} are similar to those in \eqref{e 21} on p.\ \pageref{e 21}; and again utilising \eqref{e 15}:
\begin{equation}
\Pi_P^{q}=\Pi=I-D_p\sum_{k=0}^Kq_kq_k^T\qquad\qquad\Pi^2=\Pi\qquad\qquad\Pi D_p\Pi^T
=\Pi D_p=D_p\Pi^T
\label{e 22}
\end{equation}
Thus
\[
Cov(\Delta\widehat v_P)=D_p^{1/2}\widehat\pi_{\sqrt p} D_p^{1/2}=\Pi^{q}_PD_p{\Pi^{q}_P}^T=\Pi^{q}_P D_p=D_p{\Pi^{q}_P}^T
\]
Changing to function parametric form yields
\[
\phi^T\Delta\widehat v_P
=\phi^T\left(I-D_p\sum_{k=0}^Kq_kq_k^T\right)\vec{\Delta w_P}
=\phi^T\Pi^{q}_P\vec{\Delta w_P}
=\vec{\Delta w_P}^T{\Pi^{q}_P}^T\phi
\]
Once again, as in \eqref{e 10} on p.\ \pageref{e 10}, when operating on functions in the dual space, the projection operator becomes the transpose of the operator $\Pi^{q}_P$ in the primal space.

The $q$ projected $P$ BM is
\[
v^q_P(\phi)
=\widehat v_P(\phi)
=\phi^T\Pi^{q}_P\vec{\Delta w_P}
=\phi^T\left(I-D_p\sum_{k=0}^Kq_kq_k^T\right)\vec{\Delta w_P}
\]
\[
v^q_P(\phi)=w_P(\phi)\ -\sum_{k=0}^K<\phi,q_k>_P\ w_P(q_k)
\]
which becomes, in more conventional form,
\[
\int_{-\infty}^{\infty}\phi(x)dv^q_P(x)
=\int_{-\infty}^{\infty}\phi(x)dw_P(x)
-\sum_{k=0}^K\int_{-\infty}^{\infty}\phi(x)q_k(x)dP(x)
\int_{-\infty}^{\infty} q_k(x)dw_P(x)
\]
The point parametric version is
\[
v^q_F(x)=\widehat v_P(x)=w_P(x)-\sum_{k=0}^K\int_{-\infty}^xq_k(y)dP(y)\int_{-\infty}^{\infty}q_k(x)dw_P(x)
\]

\subsubsection{Covariance of $v^q_P(\phi)$ and $v^q_P(\widetilde\phi)$ for general $q$}

\[
Cov(v^q_P(\phi),v^q_P(\widetilde\phi))
=E\,\phi^T\Pi^{q}_P\vec{\Delta w_P}\,\vec{\Delta w_P}^T\,{\Pi^{q}_P}^T\widetilde\phi
=\phi^T\Pi^{q}_PD_p{\Pi^{q}_P}^T\widetilde\phi
\]
\[
  =\phi^T\left(D_p-D_p\sum_{k=0}^Kq_kq_k^TD_p\right)\widetilde\phi
  =\phi^T\,D_p^{1/2}\,\widehat\pi_{\sqrt p}\,D_p^{1/2}\,\widetilde\phi
\]

\subsubsection{For consistency with chi squared when df = $n-K-1$}

For comparison with the mean of the chi-squared distribution when parameters were assumed known, as in \eqref{e 14} on p.\ \pageref{e 14}, we have
\[
E\ {\widehat Y}^T\,\widehat Y=E\ X^T\widehat\pi_{\sqrt p}^{\,T}\,\widehat\pi_{\sqrt p} X
=E\ X^T\widehat\pi_{\sqrt p} X
=E\left[X^TX-\sum_{k=0}^KX^TD_p^{1/2}q_kq_k^TD_p^{1/2}X\right]
\]
\[
=E\left[X^TX-\sum_{k=0}^Kq_k^TD_p^{1/2}XX^TD_p^{1/2}q_k\right]
=I-\sum_{k=0}^Kq_k^TD_pq_k
=n-(K+1)
\]

\subsection{Rotation operators in the general case}

\subsubsection{The rotation operator as a succession of reflections\label{s 3}}

We require a matrix $V_K$ with the properties that
$V_KLs_k=q_k$ for $0\leq k\leq K$, and $V_K^TD_pV_K=D_p$; i.e., the linear map sends score functions to score functions, with an adjustment to allow for the different functional spaces $L^2_P$ and $L^2_R$; and the transform also preserves the $P$-norm.

Khmaladze's method for finding a suitable $V_K$ is adumbrated in \cite{khmaladze2013a}, \cite{nguyen-ttm2017b} and \cite{kennedy2018a}, {\em i.a.}; but the methodology is set out more fully in \cite{roberts2019a}, which we follow.

\[
V_0=W_0=U_{q_0,Ls_0}\qquad \widetilde{Ls_1}=W_0Ls_1\qquad W_1=U_{q_1,\widetilde{Ls_1}}
\qquad V_1=W_1\,W_0
\]
\[
\widetilde{Ls_2}=W_1W_0Ls_2=V_1Ls_2\qquad W_2=U_{q_2,\widetilde{Ls_2}}\qquad V_2=W_2\,W_1\,W_0
\]
and so on.  More formally we have the recursion
\[
\widetilde{Ls_{j}}=V_{j-1}Ls_{j}\qquad W_j=U_{q_j,\widetilde{Ls_j}}\qquad V_j=\prod_{k=0}^jW_k\qquad \qquad \mbox{for}\qquad j\geq1
\]
Then \cite{roberts2019a} shows that $V_KLs_k=q_k$ for $0\leq k\leq K$; and it is straightforward to show that $V_K^TD_pV_K=D_p$.

\subsubsection{Rotation from one empirical process to another}

The rotation from one stochastic process to another is given by
\begin{equation}
v^{s}_R(\psi)=v^{q}_P(V_KL\psi)
\label{e 23}
\end{equation}

Proof proceeds as in \eqref{e 54} on p.\ \pageref{e 54}, again utilising \eqref{e 15}.
\[
Cov\left(v^{q}_P(V_KL\psi),v^{q}_P(V_KL\widetilde\psi)\right)=\psi^TL^TV_K^T\left(D_p-\sum_{k=0}^KD_pq_kq_k^TD_p\right)V_KL\widetilde\psi
\]
\[
=\psi^TL^TD_pL\widetilde\psi-\sum_{k=0}^K\psi^TL^TV_K^TD_pV_KLs_k\,s_k^TL^TV_K^TD_pV_KL\widetilde\psi
\]
\[
=\psi^TL^TD_pL\widetilde\psi-\sum_{k=0}^K\psi^TL^TD_pLs_k\,s_k^TL^TD_pL\widetilde\psi
\]
\[
=\psi^TD_r\widetilde\psi-\sum_{k=0}^K\psi^TD_rs_k\,s_k^TD_r\widetilde\psi
=\psi^T\left(D_r-\sum_{k=0}^KD_rs_ks_k^TD_r\right)\widetilde\psi
=Cov\left(v^{s}_R(\psi),v^{s}_R(\widetilde\psi)\right)
\]

\subsection{Higher dimensional Khmadadze rotations, and the colour-blind problem}

We consider firstly how to fit the two-dimensional Khmaladze rotation into the linear algebraic framework detailed above, and comment briefly on its three-dimensional cousin.
Largely taking our cue from
\cite{dumitrescu-khmaladze2019a}, we then symmetrise the underlying Borel sets to look at the `colour-blind' problem.

\subsubsection{The Khmadadze rotation in two dimensions}

Consider a sample of $n$ realisations of independent and identically distributed pairs $(X_j,Y_j)$, with distribution function $H(x,y)$, and marginal distribution functions $F(x), G(y)$.  We approximate $F$ and $G$ by discrete distributions $P$ and $R$ respectively, along the lines of \S \ref{s 2} on p.\ \pageref{s 2}.  This time however, we are not rotating from $F$ to $G$; rather both distributions are necessary to build the framework for our analysis.

The EDF is given by
\[
H_n(x,y)=\frac1n\sum_{j=1}^n\mathds1_{\{x-X_j\geq0\}}\mathds1_{\{y-Y_j\geq0\}}
\]
and the basic empirical process as
\[
v_n(x,y)=H_n(x,y)-H(x,y)
\]
Rather than mapping the one dimensional BM onto a BB, we now use the 2 dimensional BM $w_H(x,y)$ in `time' $H(x,y)$, projecting it (or rather tying it down at the edges) to obtain a Brownian `pillow' \citep[\S2.1]{dumitrescu-khmaladze2019a}.  Writing $w$ for $w_H$ we have:
\[
v(x,y)=w(x,y)-x\,w(\infty,y)-y\,w(x,\infty)+x\,y\,w(\infty,\infty)
\]
in which we may `anchor' the BM at the origin;  or more generally at the lower end of the supports $(M_F,M_G)$ as in \S\ref{s 2}.

The analysis of the preceding sections involving the Khmaladze rotation proceeds via
\[
v_n(\phi)=\int\phi(x,y)dv_n(x,y)\qquad\mbox{and}\qquad v(\phi)=\int\phi(x,y)dv(x,y)
\]
Now define the rectangle $R(a,b)$ in the plane as
\[
R(a,b)=\{(x,y):x\leq a,y\leq b\}
\]
Then we set $\phi_{a,b}(x,y)=\mathds1_{\{(x,y)\in R(a,b)\}}$; and the family $\Phi$ is generated by all such functions,
for $M_F<a,M_G<b$.

The central result \eqref{e 23} on p.\ \pageref{e 23} remains valid, even if the vectors of length $N$ in the previous development are now vectors of length $N^2$ as we vectorise functions and variables over the plane.

Moving to higher dimensions is straightforward, conceptually at least.
The Brownian pillow, or its equivalent, in three dimensions assumes the form
\[
z(x,y,z)=w(x,y,z)-z\,w(x,y,\infty)-y\,w(x,\infty,z)-x\,w(\infty,y,z)
\]
\[
+y\,z\,w(x,\infty,\infty)+x\,z\,w(\infty,y,\infty)+x\,y\,w(\infty,\infty,z)
-x\,y\,z\,w(\infty,\infty,\infty)
\]
Changing notation in anticipation of our work in the next section, set
\begin{equation}
R(a_1,a_2,a_3)=\{(x,y,z):x\leq a_1,y\leq a_2,z\leq a_3\}
\label{e 31}
\end{equation}
and again the family $\Phi$ is generated by all functions of the form $\phi_{a_1,a_2,a_3}(x,y,z)=\mathds1_{\{(x,y,z)\in R(a_1,a_2,a_3)\}}$.

The extension to higher dimensions is feasible in principle, but how useful the Khmaladze rotation will be in higher dimensions remains to be seen.

\subsubsection{The colour-blind problem}

In a recent article \cite{dumitrescu-khmaladze2019a} have reconsidered the colour blind problem,
in which pairs of observed coloured items are unable to be distinguished by a colour blind observer. Further background is available in 
\cite{parsadanishvili1982a} and \cite{parsadanishvili_khmaladze1982a}.

To be precise,
suppose that weights $X_r$ and $X_g$ of red and green balls have distribution functions $P_r(x)$ and $P_g(x)$ respectively, and the $i$th data point consists of a pair $\left(X_r^{(i)},X_g^{(i)}\right)$, where the random variables are independent
from data point to data point, although the weights of the red and green balls need not be independent.  Then \cite{dumitrescu-khmaladze2019a}, projecting the empirical process and utilising the properties of the Brownian pillow, consider the possibilities for statistical inference available to a colour blind person, 
able to measure the maximum and minimum weights of each pair, but without any means of knowing whether that maximum or minimum comes from the red or the green ball.

We symmetrise the rectangles $R(a,b)$ by defining
\[
S(a,b)=S(b,a)=R(a,b)\cup R(b,a)
\]
Now we set $\phi_{a,b}(x,y)=\mathds1_{\{(x,y)\in S(a,b)\}}$; the function $\phi$ is symmetric in $x$ and $y$, necessarily so by virtue of the fact that the colour blind person cannot distinguish between the colours.

The derivation of \eqref{e 23} proceeds as previously, but the vectors now have length $N(N+1)/2$.

In three dimensions, and retaining the notation in \eqref{e 31}, we have
\[
S(a_1,a_2,a_3)=\cup_{\sigma\in S_3}R(a_{\sigma1,\sigma2,\sigma3})
\]
where $S_3$ is the symmetric group on 3 symbols.  Again extension to higher dimensions, and more than two colours, is feasible conceptually; but its utility in practice, and especially in the context of the Khmaladze rotation, remains to be seen.

\section{The Khmadadze transform}

\citet[ch.\ 7]{khmaladze2013b} has derived the KT in the simple setting of a mortality investigation of $n$ independent and identically distributed lives.
After a short recapitulation of least squares regression, we shall follow his treatment closely, clarifying some points as we go.

\subsection{The linear regression model}

The model is
\[
y_t=\begin{pmatrix}x_{1t}&x_{2t}\end{pmatrix}\begin{pmatrix}a_t\\b_t\end{pmatrix}+\epsilon_t
\]
Centre this to obtain
\[
y_t-Ey_t=\begin{pmatrix}x_{1t}-Ex_{1t}&x_{2t}-Ex_{2t}\end{pmatrix}\begin{pmatrix}a_t\\b_t\end{pmatrix}+\epsilon_t
\]
from which
\[
\begin{pmatrix}x_{1t}-Ex_{1t}\\x_{2t}-Ex_{2t}\end{pmatrix}(y_t-Ey_t)
=\begin{pmatrix}x_{1t}-Ex_{1t}\\x_{2t}-Ex_{2t}\end{pmatrix}
\begin{pmatrix}x_{1t}-Ex_{1t}&x_{2t}-Ex_{2t}\end{pmatrix}
\begin{pmatrix}a_t\\b_t\end{pmatrix}
+\begin{pmatrix}x_{1t}-Ex_{1t}\\x_{2t}-Ex_{2t}\end{pmatrix}\epsilon_t
\]
Assuming the covariates and the residual are uncorrelated, taking expectations yields
\[
E\,\begin{pmatrix}x_{1t}-Ex_{1t}\\x_{2t}-Ex_{2t}\end{pmatrix}(y_t-Ey_t)
=\mbox{Covar}(x_{1t},x_{2t})\begin{pmatrix}a_t\\b_t\end{pmatrix}
\]
in which $\mbox{Covar}(x,y)$ stands for the covariance matrix, while $\mbox{Cov}(x,y)$ stands for the covariance of $x$ and $y$, the off diagonal elements of the covariance matrix.  We have now
\[
\begin{pmatrix}\widehat a_t\\\widehat b_t\end{pmatrix}
=\mbox{Covar}(x_{1t},x_{2t})^{-1}
E\,\begin{pmatrix}x_{1t}-Ex_{1t}\\x_{2t}-Ex_{2t}\end{pmatrix}(y_t-Ey_t)
\]
\[
=\mbox{Covar}(x_{1t},x_{2t})^{-1}
E\,\begin{pmatrix}x_{1t}-Ex_{1t}\\x_{2t}-Ex_{2t}\end{pmatrix}y_t
\]
Finally the predicted value of the regressand is
\begin{equation}
\widehat y_t=\begin{pmatrix}x_{1t}&x_{2t}\end{pmatrix}\mbox{Covar}(x_{1t},x_{2t})^{-1}\begin{pmatrix}E\,(x_{1t}-Ex_{1t})y_t\\E\,(x_{2t}-Ex_{2t})y_t\end{pmatrix}
\label{e 45}
\end{equation}
In econometrics, the linear predictor from regression models can seem unrelated to the last equation, although the difference is more of form than of substance.  When the parameters are not time varying, it is usual to stack the regressand and regressors, with the model expressed as $Y=X\beta+\epsilon$, with $Y$ a column vector and $X$ a matrix, whence the predicted value of $Y$ is $\widehat Y=X(X^TX)^{-1}X^TY$, with the moments being estimated from the sample.  The projection from $Y$ to $\widehat Y$ in \eqref{e 4} on p.\ \pageref{e 4}, for instance, is of this form.

The regression sought in the case of KT1 is not amenable to the stacking of variables in this way because the distributions of the random variables $x_{1t}$ and $x_{2t}$ depend on $t$, as do the coefficients.

Instead of \eqref{e 45}, the formula for the predicted value of $y_t$ is often taken to be
\begin{equation}
\widehat y_t=\begin{pmatrix}x_{1t}&x_{2t}\end{pmatrix}\mbox{Covar}(x_{1t},x_{2t})^{-1}\begin{pmatrix}E\,x_{1t}y_t\\E\,x_{2t}y_t\end{pmatrix}
\label{e 28}
\end{equation}
which is the formula used, for example in \cite{koul-swordson2011a}.

\subsection{Choice of regressand and regressors}

Following \cite[ch.\ 7]{khmaladze2013b}, we place the derivation of KT1 in the context of a mortality investigation.  There are $n$ people in the sample, and we are investigating the duration of life, so that the variable $x$ in \S\ref{s 1} is time.
The intention is to predict $\nu_l$, the number of deaths in the $l$th cell, from a linear regression model.

The empirical distribution function (EDF) is defined as $\widehat F_n(x)=\frac1n\mathds1_{\{x_j\leq x\}}$, with increment $\Delta\widehat F_n(x)=\widehat F_n(x+\Delta x)-\widehat F_n(x)$, which becomes in our previous notation $\Delta\widehat F_n(x_j)=\nu_j/n$.  Then $\widehat F(x_l)$ is the actual proportion of the sample dead by time $x_l$, and $F_\theta(x_l)=E\,\widehat F(x_l)$ the expected proportion.

Recalling the definition of the score function $Q_1$ in \eqref{e 47} on p.\ \pageref{e 47}, and assuming there is but one unknown parameter, so that $K=1$, we have
\[
\qquad Q_1=\begin{pmatrix}Q_{11}&Q_{12}&\ldots&Q_{1N}\end{pmatrix}^T
\]

To predict $\nu_l$ we would think of using its expected value $p_l$, as well as its derivative $\overset{\bullet}p_l=(\overset{\bullet}p_l/p_l)\times p_l$.  We recast these as their sample counterparts $\nu_l$ and $Q_{1l}\nu_l$,
and sum over the future cells to produce Khmaladze's regressors $\frac1n\sum_{j=l}^N\nu_j$ and $\frac1n\sum_{j=l}^N Q_{1j}(\hat\theta)\nu_j$.  The MLE $\hat\theta$ used in the second regressor is calculated from the entire sample -- see \eqref{e 30} on p.\ \pageref{e 30} below.
The scaling arises because the regressand is to be the increment in the EDF, and $\frac1n\nu_l=\Delta\widehat F_n(x_l)$.

The first regressor is an obvious enough choice, in that $\frac1n\sum_{j=l}^N\nu_j=1-\widehat F(x_l)$ is the proportion of survivors at time $x_k$; and the more surviving at time $x_l$, i.e.\ the higher the exposed to risk at time $x_l$, the greater the expected number of deaths within the period $[x_l,x_{l+1})$.

\subsection{KT1 through linear regression}

Once all the data is available, i.e.\ everyone in the sample has died, we are testing the goodness of fit of a particular distribution of duration of life, or rather a given family of distributions $F_\theta(x)$ depending on an unknown parameter $\theta$.
We assume $F_\theta(x)$ to be continuous, with
corresponding non-normalised score function denoted by $h(x,\theta)=\frac{\partial}{\partial\theta}\,f_\theta(x)\big/f_\theta(x)$.

We pretend to be partway through the cohort dying off, and try to predict the number dying in the next short interval, given firstly the number surviving at the moment, secondly the future times of death of the current survivors, and thirdly parameter estimates obtained from the entire sample. The `present' time is $t=x_l$, and we wish to predict $\nu_l$, the number dying before time $x_{l+1}$.

The key result we need is the following, taken from \citet[p.\ 79]{khmaladze2013b}
\[
\mbox{Cov}\left(\int g_1(x)\,d\widehat F_n(x),\int g_2(x)\,d\widehat F_n(x)\right)
\]
\[
=E\ \int g_1(x)\left(d\widehat F_n(x)-dF_\theta(x)\right)
\int g_2(x)\left(d\widehat F_n(x)-dF_\theta(x)\right)
\]
\begin{equation}
=\frac1n\left(\int g_1(x)g_2(x)dF_\theta(x)-\int g_1(x)dF_\theta(x)\int g_2(x)dF_\theta(x)\right)
\label{e 24}
\end{equation}
This result should be
compared with \eqref{e 44}, \eqref{e 49} and \eqref{e 25} on pp.\ \pageref{e 44}, \pageref{e 49} and \pageref{e 25} respectively; and
see also \citet[p.\ 46]{khmaladze2013b}.

In the present context,
\[
y_t=\int_t^{t+\Delta t} d\widehat F_n(x)=\int_{x_l}^{x_{l+1}} d\widehat F_n(x)\qquad\qquad
x_{1t}=\int_t^\infty d\widehat F_n(x)\qquad
\]
\[
x_{2t}=\int_t^\infty h(x,\theta) d\widehat F_n(x)\qquad\qquad
x_{2t}^*=\int_t^\infty\left[h(x,\theta)-E_\theta^t\right] d\widehat F_n(x)
\]
in which $E_\theta^t$ is the expectation of $h$ conditional upon surviving until time t:
\[
E_\theta^t=\frac{\int_t^\infty h(x,\theta)dF_\theta(x)}{1-F_\theta(t)}
\]
Expected values of $x_{2t}$ and $x_{2t}^*$ are given by
\[
E\,x_{2t}=\int_t^\infty h(x,\theta) dF_\theta(x)
=E_\theta^t[1-F_\theta(t)]
\]
and
\begin{equation}
E\,x_{2t}^*=\int_t^\infty\left[h(x,\theta)-E_\theta^t\right] dF_\theta(x)=E_\theta^t\left[1-F_\theta(t)\right]-E_\theta^t\left[1-F_\theta(t)\right]=0
\label{e 26}
\end{equation}

From \eqref{e 24} and \eqref{e 26} we have that
\[
\mbox{Cov}(x_{1t},x_{2t}^*)=\mbox{Cov}\left(\int_t^\infty d\widehat F_n(x),
\int_t^\infty\left[h(x,\theta)-E_\theta^t\right] d\widehat F_n(x)\right)
\]
\[
=\frac1n\left(\int_t^\infty\left[h(x,\theta)-E_\theta^t\right] dF_\theta(x)
-\left(1-F_\theta(t)\right)\int_t^\infty\left[h(x,\theta)-E_\theta^t\right] dF_\theta(x)\right)
\]
\[
=\frac1nF_\theta(t)\int_t^\infty\left[h(x,\theta)-E_\theta^t\right] dF_\theta(x)
=0
\]

The vanishing of $\mbox{Cov}(x_{1t},x_{2t}^*)$ simplifies calculations substantially, since the covariance matrix to be inverted in \eqref{e 45} reduces to a diagonal matrix.

Along the same lines we have
\[
\mbox{Cov}(x_{1t},x_{2t})
=\mbox{Cov}\left(\int_t^\infty d\widehat F_n(x),\int_t^\infty h(x,\theta) d\widehat F_n(x)\right)
\]
\[
=\frac1n\left(\int_t^\infty h(x,\theta)dF_\theta(x)-[1-F_\theta(t)]\int_t^\infty h(x,\theta)dF_\theta(x)\right)
\]
\[
=\frac1nF_\theta(t)\int_t^\infty h(x,\theta)dF_\theta(x)
\]
and
\[
\mbox{Cov}(y_{t},x_{2t})
=\mbox{Cov}\left(\int_t^{t+\Delta t} d\widehat F_n(x),\int_t^\infty h(x,\theta) d\widehat F_n(x)\right)
\]
\[
=\frac1n\left(\int_t^{t+\Delta t} h(x,\theta)dF_\theta(x)-[F_\theta(t+\Delta t)-F_\theta(t)]\int_t^\infty h(x,\theta)dF_\theta(x)\right)
\]
\[
\approx\frac1n\ p_l\left(h(x_l,\theta)-E_\theta^t[1-F_\theta(t)]
\right)
\]
and
\[
\mbox{Cov}(y_{t},x_{1t})
=\mbox{Cov}\left(\int_t^{t+\Delta t} d\widehat F_n(x),\int_t^\infty d\widehat F_n(x)\right)
\]
\[
=\frac1n\left(\int_t^{t+\Delta t} dF_\theta(x)-[F_\theta(t+\Delta t)-F_\theta(t)]\int_t^\infty dF_\theta(x)\right)
\]
\[
\approx\frac1n\ p_l\left(1-[1-F_\theta(t)]\right)
\]
The prediction from \eqref{e 45} on p.\ \pageref{e 45} becomes
\[
\frac1n\widehat\nu_l=\frac1n\,p_l\left[\int_t^\infty\begin{pmatrix}1&h(x_l,\theta)\end{pmatrix}d\widehat F_n(x)\right]
\mbox{Covar}(x_{1t},x_{2t})^{-1}
\begin{pmatrix}1-[1-F_\theta(t)]\\h(x_l,\theta)-E_\theta^t[1-F_\theta(t)]\end{pmatrix}
\]
The prediction from \eqref{e 28} on p.\ \pageref{e 28} becomes
\begin{equation}
\frac1n\widehat\nu_l=\frac1n\,p_l\left[\int_t^\infty\begin{pmatrix}1&h(x_l,\theta)\end{pmatrix}d\widehat F_n(x)\right]
\mbox{Covar}(x_{1t},x_{2t})^{-1}
\begin{pmatrix}1\\h(x_l,\theta)\end{pmatrix}
\label{e 29}
\end{equation}
Apart from scaling and replacing $p_l$ by $dF_\theta(t)$, the predicted value of $d\widehat F_n(x)$ as given in \eqref{e 29} agrees with that given in \cite{koul-swordson2011a} in their expression for $dw_{n\theta}(x)$, equal to $\sqrt n(y-\widehat y)$ in our notation.  Should there be further covariates, say $K$ in all, the vector $(1,h)$ would be extended to $(1,h_1,h_2,\ldots,h_K)$.

The increment $dw_{n\theta}(x)$ is in the nature of a BM, because it is uncorrelated with the past.  More precisely, $dw_{n\theta}(x)$ is uncorrelated with the future, by the nature of the regression that $\mbox{Cov}(x_{jt},y_t-\widehat y)=0$ for $j=1,2$; and the past and future are mirror images of each other.

This is clear for $x_{1t}$, since the number of survivors is the sample size minus the number who have died so far.
As for the score function,
we recall the following simple properties.
\[
\sum p_j=1\qquad\sum\overset{\bullet}p_j=0\qquad\sum\frac{\overset{\bullet}p_j}{p_j}\,p_j=0\qquad\sum Q_{1j}(\theta)p_j=0
\]
From the last of these relations, we recall that the MLE of $\theta$, say $\hat\theta$, is given by
\begin{equation}
\sum Q_{1j}(\hat\theta)\nu_j=0
\label{e 30}
\end{equation}
which is the sample counterpart of
\[
\int h(x,\theta)dF_\theta(x)=0
\]
The variable $x_{2t}$ reflects the future, but it equally reflects the past.

\section*{Acknowledgement}
Thanks to Estate Khmaladze for comments on an early draft of this paper.  Responsibility for the contents naturally remains with the author.


\end{document}